\documentclass{amsart}
\usepackage{amsmath,amsfonts,latexsym,amssymb}

\newtheorem{theorem}{Theorem}

\newtheorem{corollary}{Corollary}

\theoremstyle{remark}

\begin{document}

\title[A case of simultaneous non-vanishing]{A case of simultaneous non-vanishing}
\author{Ritabrata Munshi \and Jyoti Sengupta}
\address{School of Mathematics, Tata institute of Fundamental Research, 1 Dr. Homi Bhabha Road, Colaba, Mumbai 400005, India.} 
\email{rmunshi@math.tifr.res.in \and sengupta@math.tifr.res.in}


\begin{abstract}
We show that for $k>1000$ an even number and a sufficiently large prime $q$, there exists a newform $f$ of weight $k$ and level $q$ such that 
$$
L(1/2,f)L(1/2,\text{Sym}^2 f)\neq 0.
$$ 
\end{abstract}

\maketitle


\section{Introduction}
\label{intro}

The non-vanishing of automorphic $L$-functions at the central critical point is a topic of great interest. It has analytic, arithmetic and geometric ramifications as is clear for example from the Birch and Swinnerton-Dyer conjecture. Of no less interest is the simultaneous non-vanishing at the central critical point of two automorphic $L$-functions. This has deep implications in a wide variety of topics from Landau-Siegel zero (see \cite{IS}) to Langlands program. \\   

In this article we consider such a problem and answer it in the affirmative. Let $f$ be a primitive cusp form of weight $k$ and level $q$ (and trivial nebentypus). Let $L(s,f)$ be the degree two $L$-function associated with $f$, and let $L(s,\text{Sym}^2 f)$ be the degree three $L$-function associated with the symmetric square lift of $f$. From the work of Gelbart and Jacquet we know that this degree three $L$-function is automorphic. Let $H_k(q)$ be the collection of Hecke normalized newforms, and for $f\in H_k(q)$ we define
$$
w_f^{-1}=\frac{\Gamma(k-1)}{(4\pi)^{k-1}\|f\|^2}
$$
to be the spectral weight.
We will prove the following theorem.\\

\begin{theorem}
\label{mthm}
Let $q$ be a prime and let $k>1000$ be a large even integer. Then we have
\begin{align}
\label{cond}
\sum_{f\in H_k(q)}w_f^{-1}L(1/2,f)L(1/2,\mathrm{Sym}^2 f)=&A_k\log q+B_k+O(q^{-1/8+\varepsilon}),
\end{align}
where the constants $A_k>0$, $B_k$, as well as the implied constant, depend only on the weight $k$. (These constants are related to the gamma factors associated to $f$ and $\mathrm{Sym}^2\:f$ and are defined in the following section.)\\   
\end{theorem}

As a consequence we get the following result regarding simultaneous non-vanishing, which is in fact the main motivation for considering the above asymptotic.\\

\begin{corollary}
For $k>1000$ and a sufficiently large prime $q$, there exists a $f\in H_k(q)$ such that
\begin{align*}
L(1/2,f)L(1/2,\mathrm{Sym}^2 f)\neq 0.
\end{align*}
\end{corollary}

\bigskip

To prove the theorem we will apply the approximate functional equation to express the central values in terms of rapidly decaying Dirichlet series, and then we will apply the Petersson trace formula. The diagonal term yields the main term in the asymptotic. The main problem is to show that the off-diagonal contribution is comparatively small. To achieve this we apply the Poisson summation formula followed by evaluation of certain character sums. A satisfactory bound is then obtained by applying Heath-Brown's large sieve for quadratic characters.\\

A rather soft analysis suffices for our purpose as the problem under consideration is not actually at the `threshold of current technology'. We are studying an Euler product of degree five. To put it in proper perspective, one should note that a similar asymptotic for the degree six Euler product
\begin{align*}
\sum_{f\in H_k(q)}w_f^{-1}L(1/2,\mathrm{Sym}^2 f)^2
\end{align*}
is not yet known. This problem is just beyond the current state of technology. (Note that an asymptotic for the related sum with forms having primitive quadratic nebentypus has been established by Blomer \cite{B}.) Curiously an asymptotic for the degree eight Euler product
\begin{align*}
\sum_{f\in H_k(q)}w_f^{-1}L(1/2,f)^4
\end{align*}
is known due to deep work of Kowalski, Michel and Vanderkam \cite{KMV}. In the light of this the following problem for a degree seven Euler product seems quite interesting. Establish asymptotic for
\begin{align*}
\sum_{f\in H_k(q)}w_f^{-1}L(1/2,f)^2L(1/2,\mathrm{Sym}^2 f).
\end{align*}
Of course for the application that we noted above this is overly decorative  and hence we do not venture to tackle it in this note.

\bigskip


\section{Initial steps}
\label{is}

Let $S$ denote the sum appearing in the left hand side of \eqref{cond}.
From the approximate functional equation (see \cite{IK}) we get
\begin{align*}
L(1/2,\text{Sym}^2 f)=2\sum_{n=1}^\infty \frac{\lambda_f(n^2)}{\sqrt{n}}V\left(\frac{n}{q}\right)
\end{align*}
where 
\begin{align*}
V(y)=\frac{1}{2\pi i}\int_{(2)}\frac{\gamma_\Box(1/2+s)}{\gamma_\Box(1/2)}\:\left(1-\frac{1}{q^{1+2s}}\right)\zeta(1+2s)y^{-s}\frac{\mathrm{d}s}{s}
\end{align*}
with 
\begin{align*}
\gamma_\Box(s)=\pi^{-3s/2}\Gamma\left(\frac{s+1}{2}\right)\Gamma\left(\frac{s+k-1}{2}\right)\Gamma\left(\frac{s+k}{2}\right).
\end{align*}
The function $V$ is smooth and it satisfies
\begin{align*}
y^jV^{(j)}(y)\ll_A q^\varepsilon \left(1+|y|\right)^{-A},
\end{align*}
for any $A>1$. Also by moving the contour to the left and calculating the residue at the double pole at $s=0$ one gets
\begin{align}
\label{small-v}
V(y)=a_k\log y +b_k +O(y^{1/4}+q^{-1}),
\end{align}
for some numbers $a_k$ and $b_k$ which depend only on the weight $k$ and not on $q$. These numbers can be expressed explicitly in terms of the gamma factor $\gamma_\Box$. For our purpose we just note that $a_k\neq 0$.\\

Similarly we have the expression
\begin{align*}
L(1/2,f)=\sum_{m=1}^\infty \frac{\lambda_f(m)}{\sqrt{m}}W\left(\frac{m}{\sqrt{q}}\right)+i^k\sqrt{q}\lambda_f(q)\sum_{m=1}^\infty \frac{\lambda_f(m)}{\sqrt{m}}W\left(\frac{m}{\sqrt{q}}\right)
\end{align*}
where the function $W$ is given by
\begin{align*}
W(y)=\frac{1}{2\pi i}\int_{(2)}\frac{\gamma(1/2+s)}{\gamma(1/2)}\:y^{-s}\frac{\mathrm{d}s}{s}
\end{align*}
with 
\begin{align*}
\gamma(s)=(2\pi)^{-s}\Gamma\left(s+\frac{k-1}{2}\right).
\end{align*}
In particular $W$ is smooth, and one has 
\begin{align}
\label{w-decay}
y^jW^{(j)}(y)\ll_A q^\varepsilon \left(1+|y|\right)^{-A}.
\end{align}
Also
\begin{align}
\label{small-w}
W(y)=c_k +O(y^{1/4}),
\end{align}
for some constant $c_k\neq 0$. Note the slight difference with \eqref{small-v}, as in the previous case there was an extra pole coming from the zeta function involved in the functional equation. \\

Accordingly we have
\begin{align*}
S=2S_1+2i^k S_2
\end{align*}
where
\begin{align}
\label{s1}
S_1=\sum_{f\in H_k(q)}w_f^{-1}\sum_{n=1}^\infty \frac{\lambda_f(n^2)}{\sqrt{n}}V\left(\frac{n}{q}\right)\sum_{m=1}^\infty \frac{\lambda_f(m)}{\sqrt{m}}W\left(\frac{m}{\sqrt{q}}\right)
\end{align}
and
\begin{align}
\label{s2}
S_2=\sum_{f\in H_k(q)}w_f^{-1}\sqrt{q}\lambda_f(q)\:\sum_{n=1}^\infty \frac{\lambda_f(n^2)}{\sqrt{n}}V\left(\frac{n}{q}\right)\sum_{m=1}^\infty \frac{\lambda_f(m)}{\sqrt{m}}W\left(\frac{m}{\sqrt{q}}\right).
\end{align}
Note the mulitiplicativity relation $\lambda_f(q)\lambda_f(n^2)=\lambda_f(qn^2)$, which will be used below. Our next stop is an application of the Petersson trace formula. To this end we pick an orthogonal basis $\mathcal{B}$ of the space of old forms, and extend the above sum to a full orthogonal basis of the space $S_k(q)$. More precisely let
\begin{align*}
\mathcal{B}=\cup_{g\in H_k(1)}\{g,g^\star\}
\end{align*}
where $H_k(1)$ is the collection of normalized Hecke forms of level one, and
\begin{align*}
g^\star=g|_q-\frac{\left\langle g|_q,g\right\rangle}{\left\langle g,g \right\rangle}\:g=g|_q-\frac{\sqrt{q}\lambda_g(q)}{(q+1)}\:g.
\end{align*}
(See \S 2 of \cite{ILS} for details.)
 Consider the sum \eqref{s1} with the sum over $f$ now restricted to $\mathcal{B}$. By trivial estimation we see that this sum is bounded by $O(q^{-1/4})$. Similarly consider the modified \eqref{s2} with $f$ running over $\mathcal{B}$. A trivial estimation for this sum is not satisfactory because of the extra $\sqrt{q}$ which is coming from the root number. However
\begin{align*}
\lambda_{g^\star}(q)=\frac{1}{q^{(k-1)/2}}-\frac{\sqrt{q}\lambda_g^2(q)}{(q+1)}\ll q^{-1/2},
\end{align*}
which balances out the loss above. So the sum over $\{g^\star\}$ is satisfactory. For $f=g\in H_k(1)$ the sums over $m$ and $n$ in \eqref{s2} is bounded by $O(q^\varepsilon)$ (for example one may apply the inverse Mellin transform, move contour to the center and use the convexity bound).  \\

From the above analysis we conclude that
\begin{align*}
S_i=\Delta_i+2\pi i^{-k}\mathcal{O}_i+O(q^{-1/4}),
\end{align*}
for $i=1,2$, where $\Delta_i$ denote the diagonal contribution from the Petersson formula and $\mathcal{O}_i$ denote the off-diagonal contribution. In the rest of this section we will analyse the diagonal contributions. Note that
\begin{align*}
\Delta_1=\sum_{n=1}^\infty \frac{1}{n^{3/2}}V\left(\frac{n}{q}\right)W\left(\frac{n^2}{\sqrt{q}}\right)
\end{align*}
and
\begin{align*}
\Delta_2=\sum_{n=1}^\infty \frac{1}{n^{3/2}}V\left(\frac{n}{q}\right)W\left(\frac{qn^2}{\sqrt{q}}\right).
\end{align*}
The second diagonal is negligibly small as the function $W$ decays rapidly \eqref{w-decay}. We will now show that the first diagonal makes a sizeable contribution, and in fact it yields the main term in the asymptotic.\\

Truncating the sum at $N$ we get
\begin{align*}
\Delta_1=\sum_{n<N} \frac{1}{n^{3/2}}V\left(\frac{n}{q}\right)W\left(\frac{n^2}{\sqrt{q}}\right)+O(N^{-1/2}).
\end{align*}
Next we use \eqref{small-w} to get
\begin{align*}
c_k\sum_{n<N} \frac{1}{n^{3/2}}V\left(\frac{n}{q}\right)+O(N^{-1/2}+q^{-1/8+\varepsilon}).
\end{align*}
Then using \eqref{small-v} we obtain
\begin{align*}
a_kc_k\sum_{n<N} \frac{1}{n^{3/2}}\log (n/q)+b_kc_k\sum_{n<N} \frac{1}{n^{3/2}}+O(N^{-1/2}+q^{-1/8+\varepsilon}).
\end{align*}
Completing the remaining sums and picking $N=q^{1/4}$ we conclude that
\begin{align*}
\Delta_1=(-a_k\log q+b_k)c_k\zeta(3/2)-a_kc_k\zeta'(3/2)+O(q^{-1/8+\varepsilon}).
\end{align*}
Accordingly we set
\begin{align*}
A_k=-a_kc_k\:\zeta(3/2),\;\;\;\text{and}\;\;\;B_k=c_k(b_k\zeta(3/2)-a_k\zeta'(3/2)).
\end{align*} 

\bigskip

\section{Analysis of $\mathcal{O}_1$}
\label{o1}

To conclude the theorem it now remains to show that the off-diagonals are also bounded by $O(q^{-1/8+\varepsilon})$. We are going to establish this bound for the first off-diagonal in this section. In the next section we will show that the other off-diagonal is even smaller. The off-diagonal for the sum $S_1$ is given by
\begin{align*}
\mathcal{O}_1=\sum_{n=1}^\infty \frac{1}{\sqrt{n}}V\left(\frac{n}{q}\right)\sum_{m=1}^\infty \frac{1}{\sqrt{m}}W\left(\frac{m}{\sqrt{q}}\right)\sum_{c=1}^\infty\frac{S(m,n^2;cq)}{cq}J_{k-1}\left(\frac{4\pi n\sqrt{m}}{cq}\right).
\end{align*}
We now take a smooth dyadic subdivision of the $n$ and $m$ sums. This leads us to sums of the type
\begin{align*}
\mathcal{O}_1(N,M)=(NM)^{-1/2}\sum_{n=1}^\infty \mathcal{V}\left(\frac{n}{N}\right)\sum_{m=1}^\infty \mathcal{W}\left(\frac{m}{M}\right)\sum_{c=1}^\infty\frac{S(m,n^2;cq)}{cq}J_{k-1}\left(\frac{4\pi n\sqrt{m}}{cq}\right)
\end{align*}
with $N\ll q^{1+\varepsilon}$ and $M\ll q^{1/2+\varepsilon}$. The new functions $\mathcal{V}$ and $\mathcal{W}$ are now supported in $[1,2]$, and satisfy the bound 
\begin{align*}
\mathcal{V}^{(j)}(x),\;\mathcal{W}^{(j)}(x)\ll_j 1.
\end{align*}
Using the Weil bound for the Kloosterman sum we get that
\begin{align*}
\mathcal{O}_1(N,M)\ll q^{\varepsilon}\frac{NM^{3/4}}{q}.
\end{align*}
This is not satisfactory for our purpose.\\ 

But applying the Poisson summation formula on the $n$ sum with modulus $cq$ we will arrive at the threshold, and any more saving will be satisfactory for our purpose. Indeed after Poisson we have
\begin{align}
\label{c1c2}
\mathcal{O}_1(N,M)=N^{1/2}M^{-1/2}\sum_{n\in\mathbb{Z}} \sum_{m=1}^\infty \mathcal{W}\left(\frac{m}{M}\right)\sum_{c=1}^\infty\frac{1}{(cq)^2}\;\mathfrak{C}\;\mathfrak{I} 
\end{align}
where the character sum is given by
\begin{align*}
\mathfrak{C}=\sum_{a\bmod{cq}}S(m,a^2;cq)e\left(\frac{an}{cq}\right)
\end{align*}
and the integral is given by
\begin{align*}
\mathfrak{I}=\int_\mathbb{R} \mathcal{V}(x)J_{k-1}\left(\frac{4\pi N\sqrt{m}x}{cq}\right)e\left(-\frac{Nnx}{cq}\right)\mathrm{d}x.
\end{align*}
Extracting the oscillation of the Bessel function and then integrating by parts repeatedly we get that the integral is negligibly small if 
\begin{align*}
|n-2\sqrt{m}|\gg \frac{cq^{1+\varepsilon}}{N}.
\end{align*}\\

For simplicity let us suppose that the weight $k$ is very large (say $k>1000$), so that due to the decay of the Bessel function we are only required to consider $c$ in the range
\begin{align*}
c\ll q^\eta\frac{NM^{1/2}}{q},
\end{align*}
where $\eta>0$ is very small (for $k>1000$ we can take $\eta=1/100$). In particular we will have $(c,q)=1$. (Recall that we are also assuming that $q$ is prime.) We will now explicitly evaluate the character sum. Opening the Kloosterman sum we arrive at
\begin{align*}
\mathfrak{C}=\sideset{}{^\star}\sum_{b\bmod{cq}}\;e\left(\frac{bm}{cq}\right)\sum_{a\bmod{cq}}e\left(\frac{\bar{b}a^2+ an}{cq}\right).
\end{align*}
The inner sum is a quadratic Gauss sum, and it can be evaluated explicitly. The formula is particularly simple if $cq\equiv 1\bmod{4}$. In this case we get
\begin{align*}
\mathfrak{C}=\sqrt{cq}\sideset{}{^\star}\sum_{b\bmod{cq}}\;\left(\frac{b}{cq}\right) e\left(\frac{b(4m-n^2)}{cq}\right).
\end{align*}
Let us write $c=c_1c_2$ with $(c_1,4m-n^2)=1$ and $c_2|(4m-n^2)^\infty$. Then necessarily $c_1$ is square-free and $c_2$ is powerful, otherwise the character sum vanishes. Also we get
\begin{align}
\label{eval-char-sum}
\mathfrak{C}=\varepsilon_{qc_1}\: qc_1\sqrt{c_2}\left(\frac{4m-n^2}{c_1q}\right)\sideset{}{^\star}\sum_{b\bmod{c_2}}\;\left(\frac{b}{c_2}\right) e\left(\frac{b(4m-n^2)}{c_2}\right).
\end{align}
A similar formula can be established for other values of $c$ using the explicit formula for the generalized quadratic Gauss sum as given in Lemma~2 of \cite{B}, for example. Note that for odd $d$, we have $\varepsilon_d=1$ if $d\equiv 1\bmod{4}$ and $\varepsilon_d=i$ if $d\equiv 3\bmod{4}$. \\

Consider the sum in \eqref{c1c2} with the restriction that $c_1\sim C_1$, $c_2\sim C_2$ and 
$$
|n|\ll \mathcal{N}= q^\varepsilon(M^{1/2}+Cq/N).
$$ 
Suppose that $\mathfrak{C}$ is replaced by the explicit expression \eqref{eval-char-sum}. The resulting sum will be denoted by $\mathcal{O}_1(N,M,C_1,C_2)$. We will now obtain a sufficient bound for this sum, see \eqref{bound-for-o1} below, which will be uniform with respect to all the parameters. It will be clear that this bound holds for $\mathcal{O}_1$ as well. Let $C=C_1C_2$. By a change of variables we have
\begin{align*}
\mathfrak{I}=\frac{c}{C}\:\int_\mathbb{R} \mathcal{V}\left(\frac{cx}{C}\right) J_{k-1}\left(\frac{4\pi N\sqrt{m}x}{Cq}\right)e\left(-\frac{Nnx}{Cq}\right)\mathrm{d}x.
\end{align*} 
This separates the variable $c$ from $(m,n)$. Taking absolute values we get
\begin{align*}
\mathcal{O}_1(N,M,& C_1,C_2)\ll \frac{q^\varepsilon N^{1/2}}{qM^{1/2}C_1C_2^{1/2}}\sum_{n\ll \mathcal{N}}\; \sum_{m\sim M}\\
&\times \sum_{\substack{c_2\sim C_2\\c_2|(4m-n^2)^\infty\\c_2\:\text{powerful}}}\min\left\{\left(\frac{NM^{1/2}}{Cq}\right)^{11},\frac{\sqrt{Cq}}{\sqrt{N}M^{1/4}}\right\}\left|\sum_{\substack{c_1\sim C_1\\c_1\:\Box-\text{free}}}\;\varepsilon_{qc_1}\:\left(\frac{4m-n^2}{c_1}\right)\right|
\end{align*}
where $C_1C_2=C\ll q^{\eta-1}NM^{1/2}$. (Note that the factor before the absolute value sign takes into account the size of the Bessel function.) Consider the terms where $4m=n^2$. In this case $C_1=1$ and the innermost sum is trivial. Evaluating the other sums trivially we see that this makes a contribution of size $O(q^{\eta-1/2})$. Hence
\begin{align*}
\mathcal{O}_1(N,M,C_1,C_2)\ll &q^\varepsilon\min\left\{\frac{N^{23/2}C_2^{1/2}}{q^{12}C^{12}},\frac{1}{C_1^{1/2}q^{1/2}M^{3/4}}\right\}\sup_U\sum_{u\sim U}\\
&\times \sum_{\substack{1\leq d\ll D\\d\:\Box-\text{free}}}\alpha(d,u)\left|\sum_{\substack{c_1\sim C_1\\c_1\:\Box-\text{free}\\(c_1,u)=1}}\;\varepsilon_{qc_1}\:\left(\frac{d}{c_1}\right)\right|+q^{\eta-1/2}
\end{align*}
where $D=Cq^{1+\varepsilon}\mathcal{N}/NU^2$, and
\begin{align*}
\alpha(d,u)=\mathop{\sum_{n\ll \mathcal{N}}\: \sum_{m\sim M}}_{4m-n^2=du^2}\;\sum_{\substack{c_2\sim C_2\\c_2|(4m-n^2)^\infty\\c_2\:\text{powerful}}}\;1\ll \mathcal{N} q^{\varepsilon}.
\end{align*}\\

By Cauchy inequality we arrive at the following expression
\begin{align*}
\sum_{\substack{1\leq d\ll D\\d\:\Box-\text{free}}}\left|\sum_{\substack{c_1\sim C_1\\c_1\:\Box-\text{free}\\(c_1,u)=1}}\;\varepsilon_{qc_1}\:\left(\frac{d}{c_1}\right)\right|^2.
\end{align*}
Using Heath-Brown's large sieve inequality (see \cite{H}) we dominate the above expression by
\begin{align*}
q^\varepsilon\left(D+C_1\right)C_1.
\end{align*}
Consequently it follows that
\begin{align*}
\mathcal{O}_1(N,M,C_1,C_2)\ll q^{\eta-1/2}+q^\varepsilon &\min\left\{\frac{N^{23/2}C_2^{1/2}}{q^{12}C^{12}},\frac{1}{C_1^{1/2}q^{1/2}M^{3/4}}\right\}\mathcal{N}\\
&\times \sup_U\sum_{u\sim U}\:D^{1/2}\:\left(D+C_1\right)^{1/2}C_1^{1/2}.
\end{align*}
Hence
\begin{align}
\label{bound-for-o1}
\mathcal{O}_1(N,M,C_1,C_2)\ll q^{\eta-1/2}+q^{-1/8+\varepsilon}\ll q^{-1/8+\varepsilon}.
\end{align}
Since the last bound holds for all values of $N$ and $M$, we get that
\begin{align}
\label{bound-o1}
\mathcal{O}_1\ll q^{-1/8+\varepsilon}.
\end{align}\\


\section{Analysis of $\mathcal{O}_2$}
\label{o2}

In this section we will analyse the second off-diagonal.
The off-diagonal for the sum $S_2$ is given by
\begin{align*}
\mathcal{O}_2=\sqrt{q}\sum_{n=1}^\infty \frac{1}{\sqrt{n}}V\left(\frac{n}{q}\right)\sum_{m=1}^\infty \frac{1}{\sqrt{m}}W\left(\frac{m}{\sqrt{q}}\right)\sum_{c=1}^\infty\frac{S(m,qn^2;cq)}{cq}J_{k-1}\left(\frac{4\pi n\sqrt{m}}{c\sqrt{q}}\right).
\end{align*}
As before we take a smooth dyadic subdivision of the $n$ and $m$ sums. This leads us to sums of the type
\begin{align*}
\mathcal{O}_2(N,M)=(NM)^{-1/2}\sum_{n=1}^\infty \mathcal{V}\left(\frac{n}{N}\right)\sum_{m=1}^\infty \mathcal{W}\left(\frac{m}{M}\right)\sum_{c=1}^\infty\frac{S(m,qn^2;cq)}{cq}J_{k-1}\left(\frac{4\pi n\sqrt{m}}{c\sqrt{q}}\right).
\end{align*}
Here $N\ll q^{1+\varepsilon}$, $M\ll q^{1/2+\varepsilon}$ and the functions $\mathcal{V}$, $\mathcal{W}$ are as in the previous section. Here the $c$ sum can be truncated at $c\ll q^{1/100-1/2}NM^{1/2}$ at a cost of an error term of size $O(q^{-1})$.\\ 

We now apply the Poisson summation formula on the $n$ sum with modulus $c$. (Observe the drop in the modulus due to form of the root number.) This yields
\begin{align*}
&\sum_{n=1}^\infty S(m,qn^2;cq)J_{k-1}\left(\frac{4\pi n\sqrt{m}}{c\sqrt{q}}\right)\mathcal{V}\left(\frac{n}{N}\right)=\frac{N}{c}\sum_{n\in\mathbb{Z}}\;\mathfrak{C}\:\mathfrak{I}
\end{align*}
where the new character sum is given by
\begin{align*}
\mathfrak{C}=\sum_{a\bmod{c}}S(m,qa^2;cq)e\left(\frac{an}{c}\right)
\end{align*}
and the integral is given by
\begin{align*}
\mathfrak{I}=\int_\mathbb{R} \mathcal{V}\left(x\right)J_{k-1}\left(\frac{4\pi N\sqrt{m}x}{c\sqrt{q}}\right)e\left(-\frac{Nnx}{c}\right)\mathrm{d}x.
\end{align*}
Using the decomposition of the Bessel function and integrating by parts we get that the integral is negligibly small if 
\begin{align*}
|\sqrt{m}-\sqrt{q}n|\gg q^\varepsilon\frac{c}{N}.
\end{align*}
If $n\neq 0$ then the left hand side is at least of the size $q^{1/2}$. But if $c\gg Nq^{1/4}$ then $N\sqrt{m}/c\sqrt{q}\ll q^{-1/2+\varepsilon}$. Consequently 
\begin{align*}
\mathfrak{I}\ll q^{-k/4},
\end{align*} 
and for $k$ large enough this gives a satisfactory bound even after estimating the other sums trivially. This takes care of the contribution of $n\neq 0$. For $n=0$ the inequality $\sqrt{m}\ll q^\varepsilon c/N$ again implies that $N\sqrt{m}/c\sqrt{q}\ll q^{-1/2+\varepsilon}$, and consequently the same bound holds in this case as well. We conclude that 
\begin{align}
\label{bound-o2}
\mathcal{O}_2\ll q^{-1/8}.
\end{align} 
Actually we can conclude a much stronger bound as we are picking $k>1000$. But this is sufficient for the present purpose.

\bigskip


\end{document}